\documentclass[10pt]{IEEEtran}

\ifCLASSINFOpdf
\else
   \usepackage[dvips]{graphicx}
\fi
\usepackage{url}

\usepackage{graphicx}
\usepackage[ruled,vlined]{algorithm2e}
\usepackage{amsfonts,amsthm,amsmath,xcolor,amssymb}
\usepackage{diagbox}
\usepackage{ulem}
\usepackage{algorithmic}

\newtheorem{prop}{Proposition}

\theoremstyle{definition}

\newtheorem{rem}{Remark}
\newtheorem{prob}{Problem}

\usepackage{flushend}

\newcommand{\prior}{\hat{\pi}}
\newcommand{\du}{\delta_U}
\newcommand{\dl}{\delta_L}
\newcommand{\gu}{\gamma_U}
\newcommand{\gl}{\gamma_L}
\newcommand{\bl}{\bar{\lambda}}

\newcommand{\bs}{\mathbf{s}}

\newcommand{\mylam}{\Lambda}

\newcommand{\mystop}{\hat{t}}

\begin{document}

\title{Anomaly Search Over Many Sequences \\ With Switching Costs}

\author{Matthew Ubl, \IEEEmembership{Student Member, IEEE}, Benjamin D. Robinson, \IEEEmembership{Member, IEEE}, \\
and Matthew Hale, \IEEEmembership{Member, IEEE}
\thanks{Date submitted: 2 March 2023. AFRL Public Affair number AFRL-2023-0996. This work was supported by 
the Office of Naval Research under grants
N00014-22-1-2435 and N00014-21-1-2495, by
the Air Force Office of Scientific Research under grants FA9550-19-1-0169 and RYCOR036, and by the Air Force Research Lab Sensors Directorate.}
\thanks{Benjamin Robinson is with the Air Force Research Lab Sensors Directorate, WPAFB, OH 45433 USA.  (email: benjamin.robinson.8@us.af.mil)}
\thanks{Matthew Ubl and Matthew Hale are with the University of Florida Aerospace Engineering Department, Gainesville, FL, 32611 USA.  (Respective emails: m.ubl@ufl.edu and matthewhale@ufl.edu)}}

\maketitle

\begin{abstract}
This paper considers the quickest search problem to identify anomalies among large numbers of data streams. These streams can model, for example, disjoint regions monitored by a mobile robot. A particular challenge is a version of the problem in which the experimenter must suffer a cost each time the data stream being sampled changes, such as the time the robot must spend moving between regions. In this paper, we propose an algorithm which accounts for switching costs by varying a confidence threshold that governs when the algorithm switches to a new data stream. 
Our main contributions are easily computable approximations for both the optimal value of
this threshold and the optimal value of the parameter that determines when a stream
must be re-sampled. Further, we  empirically 
show (i) a uniform improvement for switching costs of interest 
and (ii) roughly equivalent performance for small switching  costs when
comparing to the closest available algorithm.
\end{abstract}


\IEEEpeerreviewmaketitle

\section{Introduction}

\IEEEPARstart{A}{}fundamental problem in sensing and signal processing is online anomaly search.  With origins in Chernoff's sequential design of experiments \cite{chernoff1959sequential}, the aim of online anomaly search is to develop an efficient policy for sampling a subset of several data streams over time that quickly and accurately identifies any anomalous ones.  Some variations on Chernoff's method include \cite{dragalin1996simple,nitinawarat2013controlled,naghshvar2013active,nitinawarat2015controlled,cohen2015active,huang2018active,tsopelakos2019sequential}. An important application is clinical-trial design, where the drug tested in each trial depends on the drugs tested and outcomes obtained in all previous trials of the study.  
Another application is online search for an open channel in cognitive radio, where the sampling policy may depend on past observations.  


Recently several authors have considered the problem of online anomaly search in which a cost is incurred any time the current channel sampled differs from the last one \cite{vaidhiyan2017neural, vaidhiyan2015active, lambez2021anomaly}.  This new formulation models online anomaly search more realistically than the traditional formulation because it accounts for possible switching costs in hardware or software, e.g., if an autonomous agent must move to observe a new location or if equipment must be repositioned to do so. This formulation is also more challenging because conventional analyses of exploration versus exploitation 
do not apply.
Table~\ref{tab:context} shows how the problem setup in this paper relates to that of existing work.
We emphasize that we differ from these existing works by considering switching costs and an infinite number of data streams simultaneously.  To the best of our knowledge, this is the first work to do so.

Our approach is to adapt an existing method to the case of switching costs.  We do this by (i) introducing a new  parameter that controls switching between data streams and (ii) numerically optimizing both this switching parameter and an existing stopping parameter.
Our contributions are:
\begin{itemize}
\item We develop an algorithm that generalizes the algorithm in \cite{lai2011quickest}, which is known to be Bayes-optimal for the setting without switching costs.
\item We show that in a certain asymptotic regime, the optimal parameters for this algorithm can be found by exactly solving an algebraic equation and numerically solving a strongly convex optimization problem over a scalar decision variable.
This approach compares favorably with the Monte Carlo method of \cite{lai2011quickest} in terms of computational efficiency.
\item We illustrate an almost uniform improvement in performance over the closest existing method in experiments. We illustrate empirically a uniform improvement for switching costs of interest and roughly equivalent performance for small switching costs, compared to \cite{lai2011quickest}.
\end{itemize}

\begin{figure}
    \centering
\begin{tabular}{|c||c|c|} \hline
     &  Switching costs  & No switching costs \\ \hline\hline
  Finite \# seqs   &  \cite{lambez2021anomaly}   & \cite{chernoff1959sequential} \cite{dragalin1996simple} \\ \hline
  Infinite \# seqs     &  This work    & \cite{lai2011quickest} \\ \hline
\end{tabular}
    \caption{Relationship of this work (lower left quadrant) to existing works; the other
    three quadrants contain references to representative existing works
    on those classes of problems.}
    \label{tab:context}
\end{figure}

The rest of this paper is organized as follows. In Section~\ref{sec:prelim} we provide our problem setup and discuss the existing optimal solution for the setting without switching costs. In Section~\ref{sec:probstate} we introduce switching costs and derive an approximation of the combined observation-switching cost, and derive the parameter choices that minimize this approximation. In Section~\ref{sec:discussion} we justify this approximation and explore the conditions for which it is accurate. In Section~\ref{sec:simulations} we provide numerical results, and we conclude in Section~\ref{sec:con}.

\section{Preliminaries} \label{sec:prelim}
\subsection{Problem Setup}
We consider data streams indexed over~$k \in \{1, 2,  \ldots\}$. 
Data stream $k$ generates the i.i.d. sequence of random variables $X^k_1, X^k_2,\dots$ which have sample space $\mathcal{X}$. Each data stream obeys one of two hypotheses, $H_0$ or $H_1$. Consider two distinct distributions on $\mathcal{X}$: $F_0$ and $F_1$. We say hypothesis $H_0$ is true for data stream $k$ if $X^k_t \sim F_0$ for all $t \in \{1,2,\dots\}$, and that $H_1$ is true for data stream $k$ if $X^k_t \sim F_1$ for all $t \in \{1,2,\dots\}$. We use $f_0$ and $f_1$ to denote the PDFs of $F_0$ and $F_1$ respectively. In this paper if $H_0$ is true for a particular data stream, we say that stream is \textit{nominal}, and if $H_1$ is true, that the stream is a \textit{target} stream. We assume that for any given data stream, $H_1$ is true with prior probability $\prior$ and $H_0$ is true with prior probability $1-\prior$, where $\prior \in (0,1)$. We use $\mathbb{E}_{i}$ to denote the expectation under hypothesis $i$, and define $\mathbb{E}_{\prior}[\,\cdot\,] = \prior \mathbb{E}_1[\,\cdot\,] + (1-\prior) \mathbb{E}_0[\,\cdot\,]$.

We assume we have a single observer that can sample one and only one data stream at a time. That is, if the observer samples data stream $k$ at time $t$ it receives $X^k_t$ but no information from the other data streams. Our goal is to design an algorithm for this observer to identify a target data stream as quickly as possible. In much of the existing literature, ``as quickly as possible'' means minimizing the expected number of observations required while satisfying some constraint on the error probability, i.e., 
satisfying an upper bound on the probability
that the stream we identify as a target is actually nominal. We use $\tau$ to denote the number of observations taken before the algorithm terminates and declares a particular data stream, which we denote as $k_{\tau}$, as a target. We use $H^{k_{\tau}}$ to denote the hypothesis obeyed by data stream $k_{\tau}$, and $P(H^{k_{\tau}} = H_0)$ as the error rate, i.e., the probability that the stream $k_\tau$ is actually nominal. Therefore, our goal is to find the algorithm that minimizes $\mathbb{E}_{\prior}[\tau]$ while ensuring $P(H^{k_{\tau}} = H_0) \leq \epsilon$, where $\epsilon > 0$ is an allowable error rate. 


\subsection{Solution Without Switching Costs}
This subsection briefly reviews related work on problems without switching costs; switching costs will be introduced in the next section. It was shown in~\cite{lai2011quickest} that a cumulative-sum-based (CUSUM-based) test is the optimal algorithm for the setting without switching costs. This algorithm is defined by two threshold parameters: $\gl \leq 0 \leq \gu$. In this algorithm the observer maintains a statistic $\mylam^k_t$ for stream $k$ which is updated after every observation and is initialized as $\mylam^1_0 = 0$. If at time $t$ the observer samples data stream $k$ (i.e., the observer receives $X^k_t$), then it performs the update $\mylam^k_t = \mylam^k_{t-1} + \log \left(\frac{f_1(X^k_t)}{f_0(X^k_t)}\right)$. If $\gl \leq \mylam^k_t < \gu$, then the observer will sample data stream $k$ again at time $t+1$. If $\mylam^k_t < \gl$, then the observer has declared stream $k$ as nominal. The observer will begin using the new statistic $\mylam^{k+1}_t=0$ and will switch to sampling data stream $k+1$ (we assume the streams are either pre-ordered or the next stream is selected at random) beginning at time $t+1$. This procedure describes the $k^{th}$ \textit{stage} of the algorithm, during which data stream $k$ is observed. 

The algorithm carries out the same procedure on data streams $k+1,k+2,\dots$ until a target stream is indicated. Specifically, if $\mylam^k_t \geq \gu$, then the algorithm terminates and the observer declares data stream $k$ as a target stream. The optimal choice for $\gl$ is $0$ regardless of $F_0$, $F_1$, $\prior$, or $\epsilon$~\cite[Section IV]{lai2011quickest}. Because $\mathbb{E}_{\prior}[\tau]$ monotonically increases as $\gu \rightarrow \infty$ and $P(H^{k_{\tau}} = H_0)$ monotonically decreases as $\gu \rightarrow \infty$, the optimal choice for $\gu$ is the smallest value for which $P(H^{k_{\tau}} = H_0) \le \epsilon$. However, a closed form for this $\gu$ is not known; it must be estimated using numerical experiments.

\section{Problem Statement with Switching Costs} \label{sec:probstate}
We now introduce switching costs to the model described in the previous section, and this will give the problem formulation that we consider
in this paper. 
When the observer switches from data stream $k$ to data stream $k+1$, it now incurs a cost $\lambda_k$ drawn from some non-negative distribution $L$ for all $k$. That is, $\lambda_k \geq 0$ and $\mathbb{E}[\lambda_k] = \bl$ is finite. 
 We also assume that the costs $\lambda_k$ and the observations $X^{k'}_t$ are mutually independent for all $k$, $k'$, and $t$. This cost models applications in which observing a new data stream requires ``deadtime'' when no observations can be taken, such as when equipment needs to be re-positioned or re-calibrated. It also models problems in which observations and switches are ``costly'' in some resource other than time, such as energy or money. Under this switching cost assumption, the new optimization problem becomes:
\begin{prob} \label{prob1}
Let an error tolerance $\epsilon \in (0,1)$ and prior $\prior \in (0,1)$ be given. Then
\begin{align}
    &\text{minimize}_{\gl \leq 0 \leq \gu} \mathbb{E}_{\prior} [\tau] + \mathbb{E}_{\prior}[s] \label{eqn:optprob}\\
    &\text{s.t. } P(H^{k_{\tau}} = H_0) \leq \epsilon, \label{eqn:constraint}
\end{align}
\noindent where $s = \sum_{i=1}^{k_{\tau}-1} \lambda_k$ is the total switching cost incurred before terminating the algorithm.
\end{prob}

To derive threshold choices for this problem we will next rewrite the problem in terms of the \textit{stage-wise} false-positive and false-negative rates $\alpha$ and $\beta$ of the algorithm, and then establish the relationship between these rates and the threshold choices $\gl$ and $\gu$.  
By the stage-wise false-positive rate, we mean the probability that a stage's 
terminal value of $\mylam$ exceeds (or equals) $\gu$ given that the stage's data follow $H_0$, i.e., the probability that a stream is declared a target even though the stage's data are nominal.
The stage-wise false-negative rate is similarly defined as the probability that $\gl$ exceeds the stage's terminal value of $\mylam$ given that the stage's data follow $H_1$, i.e., the stage's data are not nominal, but the stream is labelled nominal. 

Let $t_k$ be the time when the algorithm takes its last observation of stream $k$ (i.e., $\mylam^k_{t_k} \notin [\gamma_L, \gamma_U))$. 
Then $\hat{t} = t_k-t_{k-1}$ is the \textit{stage-wise stopping time} of stage $k$, or the number of observations taken of stream $k$ before making a decision. Then we may more compactly write that $\alpha = \mathbb{P}_0[\Lambda^k_{t_k} \geq \gu]$ and $\beta = \mathbb{P}_1[\Lambda^k_{t_k} < \gl]$, where $\mathbb{P}_i$ is the probability under hypothesis $i$.  
Consider the inequalities 
\begin{align}
   & \gamma_L \ge \log (\beta/(1-\alpha)) \label{eq:gl-approx} \\
   & \gamma_U \le \log((1-\beta)/\alpha)) \label{eq:gu-approx} \\
   & \mathbb{E}_{\prior}[\Lambda^k_{t_k} \mid \Lambda^k_{t_k} \ge \gu] \ge \gu \label{eq:gu-ineq} \\
& \mathbb{E}_{\prior}[\Lambda^k_{t_k} \mid \Lambda^k_{t_k} < \gl] \le \gl, \label{eq:gl-ineq}
\end{align}
which are given in \cite[Equations~(2.9) and (2.10)]{siegmund1985sequential}. In accordance with \cite{siegmund1985sequential}, we assume these inequalites are approximate equalities for the remainder of this section.  These approximations are known as ``Brownian motion approximations''.
\begin{rem}[Brownian Motion Approximations]
The assumption that \eqref{eq:gl-approx}-\eqref{eq:gl-ineq} are approximate equalities is accurate under the conditions that (a) $f_0$ and $f_1$ are sufficiently close, (b) $\bl$ is sufficiently large, and (c) $\epsilon$ is small.  These conditions imply that the step size in the sequential probability ratio test of one stage is small compared to the decision thresholds, and thus overshoots of decision thresholds are also comparatively small.  We will elaborate on these conditions in Section~\ref{sec:discussion}. 
\end{rem}

The quantities $\alpha$ and $\beta$ appear in Problem~1 in the following manner: using Wald's Identity 
we see that 
$\mathbb{E}_{\prior} [\tau] = \mathbb{E}_{\prior}[k_{\tau}] \mathbb{E}_{\prior}[\mystop]$. 
Intuitively, this result states that the expected number of total observations before termination is equal to the expected number of data streams visited ($k_{\tau}$) multiplied by the expected number of observations per data stream (which is $\mystop$). Using Wald's Identity again gives $\mathbb{E}_{\prior} [s] = \mathbb{E}_{\prior}[k_{\tau}-1] \bl$, which states that the expected total switching cost incurred over time is equal to the expected number of switches (one less than the number of streams visited) multiplied by the expected switching cost per switch.

We first address the term $\mathbb{E}_{\prior}[\mystop]$. Assume that a particular data stream $k$ being sampled by the observer is a target, and assume that the observer takes its last sample of $k$ at time $t_k$. Then one of the termination criteria has been met and~$\Lambda^k_{t_k} \not\in [\gl, \gu)$. 
Let $f_0,f_1$ be the PDFs of $F_0,F_1$ respectively. 
Using Wald's Identity again, we have $\mathbb{E}_{1} [\Lambda^k_{t_k}] = D(f_1||f_0) \mathbb{E}_1[t_k]$, where $D(f_1||f_0) = \mathbb{E}_1\left[\log{\frac{f_1(x)}{f_0(x)}}\right]$ is the Kullback-Leibler (KL) divergence of $f_0$ from $f_1$. Using the same procedure we see $\mathbb{E}_{0} [\Lambda^k_{t_k}] = -D(f_0||f_1) \mathbb{E}_0[\mystop]$. From the definition of $\mathbb{E}_{\prior}[\cdot]$, we can write
\begin{equation}
\mathbb{E}_{\prior}[\mystop] = (1-\prior)\frac{\mathbb{E}_{0} [\Lambda^k_{t_k}]}{-D(f_0||f_1)} + \prior \frac{\mathbb{E}_{1} [\Lambda^k_{t_k}]}{D(f_1||f_0)}.
\end{equation}
Furthermore, we see that $\mathbb{E}_{0} [\Lambda^k_{t_k}] = \alpha \mathbb{E}_{0} [\Lambda^k_{t_k} | \Lambda^k_{t_k} \geq \gu] + (1-\alpha) \mathbb{E}_{0} [\Lambda^k_{t_k} | \Lambda^k_{t_k} < \gl]$. By the Brownian motion approximations we have $\mathbb{E}_{0} [\Lambda^k_{t_k} | \Lambda^k_{t_k} \geq \gu]  \approx \gu$ and $\mathbb{E}_{0} [\Lambda^k_{t_k} | \Lambda^k_{t_k} < \gl]  \approx \gl$.
Following equivalent steps for $\mathbb{E}_{1} [\Lambda^k_{t_k}]$ gives
\begin{equation} \label{eqn:alphagamma}
\mathbb{E}_{\prior}[\mystop] \approx (1-\prior)\frac{\alpha \gu + (1-\alpha)\gl}{-D(f_0||f_1)} + \prior \frac{(1-\beta)\gu + \beta \gl}{D(f_1||f_0)}.
\end{equation}
Furthermore, from~\cite[Equation (30)]{lai2011quickest} we also have 
\begin{equation} \label{eqn:switchesalphabeta}
\mathbb{E}_{\prior} [k_{\tau}] = \frac{1}{(1-\prior)\alpha + \prior(1-\beta)}.
\end{equation}

For ease of notation, we now define $\du = \exp (\gu)$ and $\dl = \exp (\gl)$. While $\gl \leq 0 \leq \gu$ are the actual thresholds used by the algorithm, rewriting the problem in terms of $0 < \dl \leq 1 \leq \du$ makes the following notation simpler. 
Inverting the approximate equalities in \eqref{eq:gl-approx} and \eqref{eq:gu-approx}, we obtain the following:
$\alpha \approx \frac{1-\dl}{\du-\dl}$, $\beta \approx \dl\frac{\du-1}{\du-\dl}$, $1-\alpha \approx \frac{\du-1}{\du-\dl}$, and $1-\beta \approx \du\frac{1-\dl}{\du-\dl}$. 
Substituting these into~\eqref{eqn:alphagamma} and~\eqref{eqn:switchesalphabeta} and simplifying 
yields that the function
\begin{multline} 
  C(\dl,\du) := \frac{1-\prior}{-D(f_0||f_1)} \frac{\log(\du)+\frac{\du-1}{1-\dl}\log(\dl)}{1+\prior(\du-1)} \\ \qquad\qquad\qquad\quad + \frac{\prior}{D(f_1||f_0)} \frac{\du\log(\du)+\dl\frac{\du-1}{1-\dl}\log(\dl)}{1+\prior(\du-1)} \\ 
  + \frac{\bl\frac{\du-\dl}{1-\dl}}{1+\prior(\du-1)} \label{eqn:cost}
\end{multline}
approximates the cost given in \eqref{eqn:optprob}.


Furthermore, from~\cite[Equation~(29)]{lai2011quickest}  we have that that $P(H^{k_{\tau}} = H_0) = \frac{(1-\prior)\alpha}{(1-\prior)\alpha + \prior(1-\beta)}$, which is bounded above by $ \frac{1-\prior}{1+\prior(\du-1)}$  since $(1-\beta)/\alpha \ge \delta_U$. As a result, this inequality is an approximate equality under the Brownian motion approximations. Therefore, following some algebraic manipulation we see that $P(H^{k_{\tau}} = H_0) \leq \epsilon$ if $\du \geq \frac{1-\prior}{\prior}\frac{1-\epsilon}{\epsilon}$. Note that $\frac{1-\prior}{\prior}\frac{1-\epsilon}{\epsilon} > 1$ so long as $\epsilon < 1-\prior$, which should always be the case; if the tolerable error rate is greater than the prevalence of nominal data streams then the optimal algorithm is to take no observations and flag a data stream at random.

Therefore, we can find an approximate solution to Problem~\ref{prob1} by solving the following problem.
\begin{prob}
Let an error tolerance $\epsilon \in (0,1)$ and a prior $\prior \in (0,1)$ be given. Then find
\begin{align}
 \text{Find } (\dl^*,\du^*) = \underset{\dl,\du}{\arg\min }  \text{ } & C(\dl,\du) \label{approxprob}\\
 \text{s.t. } &\du \geq \frac{1-\prior}{\prior}\frac{1-\epsilon}{\epsilon} \\
& \dl \in [0,1].
\end{align}
\end{prob}

From the structure of $C(\dl,\du)$, we can derive the following propositions:
\begin{prop} \label{prop1}
For any fixed $\hat{\delta}_L \in (0,1)$, $C(\hat{\delta}_L,\cdot)$ is monotonically increasing on the domain $[1,\infty)$. From this fact, we get $\du^* = \frac{1-\prior}{\prior}\frac{1-\epsilon}{\epsilon}$.
\end{prop}
\textit{Proof:} The monotonic behavior of $C(\hat{\delta}_L,\cdot)$ on this domain is apparent by inspection. Because we wish to minimize $C(\hat{\delta}_L,\cdot)$, we want to set $\du$ to its minimum allowable value. Since that value is $\frac{1-\prior}{\prior}\frac{1-\epsilon}{\epsilon}$, we have $\du^* = \frac{1-\prior}{\prior}\frac{1-\epsilon}{\epsilon}$ regardless of the value of $\hat{\delta}_L$. $\hfill\blacksquare$
\begin{prop} \label{prop2}
    For any fixed $\hat{\delta}_U > 1$, $C(\cdot,\hat{\delta}_U)$ is strongly convex on the domain $(0,1]$. From this fact $\dl^*$ exists, is unique, and can be found by solving the scalar optimization problem $\dl^* = \underset{\dl \in [0,1]}{\arg \min} \text{ } C(\dl,\du^*)$. Furthermore, $\dl^*$ lies in the interior of this interval for $\bl > 0$ (i.e., $\dl^* \in (0,1))$. 
\end{prop}


\textit{Proof:} Differentiation shows $\lim_{\dl \rightarrow 0^+} \frac{\partial C(\dl,\hat{\delta}_U)}{\partial \dl} = -\infty$ and $\lim_{\dl \rightarrow 1^-} \frac{\partial C(\dl,\hat{\delta}_U)}{\partial \dl} = \infty$ so long as $\bl > 0$. Therefore, from the Intermediate Value Theorem there must exist some value $\hat{\delta}_L \in (0,1)$ for which $\frac{\partial C(\dl,\hat{\delta}_U)}{\partial \dl} = 0$ at $\dl = \hat{\delta}_L$. Strong convexity is established by characterizing the limiting behavior of $\frac{\partial^2 C(\dl,\hat{\delta}_U)}{\partial \dl^2}$. Observe that $C(\dl,\du)$ is a sum of three terms: the first two which contain the KL divergences $D(f_0||f_1)$ and $D(f_1||f_0)$, and the third which contains $\bl$. Name these terms $C_1$ $C_2$, and $C_3$ respectively. 
First, we see that $\lim_{\dl \rightarrow 0^+} \frac{\partial^2 (C_1(\dl,\hat{\delta}_U) + C_2(\dl,\hat{\delta}_U))}{\partial \dl^2} = \infty$ and $\lim_{\dl \rightarrow 1^-} \frac{\partial^2 (C_1(\dl,\hat{\delta}_U) + C_2(\dl,\hat{\delta}_U))}{\partial \dl^2} = \frac{\hat{\delta}_U-1}{1+\prior(\hat{\delta}_U-1)}(\frac{2}{3}\frac{1-\prior}{D(f_0||f_1)}+\frac{1}{3}\frac{\prior}{D(f_1||f_0)})$, and that $\frac{\partial^2 (C_1(\dl,\hat{\delta}_U) + C_2(\dl,\hat{\delta}_U))}{\partial \dl^2}$ is monotonically decreasing with $\dl$ on~$[0, 1]$.
Additionally, we see $\lim_{\dl \rightarrow 0^+} \frac{\partial^2 C_3(\dl,\hat{\delta}_U)}{\partial \dl^2} = \frac{\bl(\hat{\delta}_U-1)}{1+\prior(\hat{\delta}_U-1)}$ and $\lim_{\dl \rightarrow 1^-} \frac{\partial^2 C_3(\dl,\hat{\delta}_U)}{\partial \dl^2} = \infty$, and that $\frac{\partial^2 C_3(\hat{\delta}_L,\hat{\delta}_U)}{\partial \dl^2}$ is monotonically increasing with $\dl$ on~$[0,1]$.
Therefore, $C(\,\cdot\,, \hat{\delta}_U)$ is $\frac{\hat{\delta}_U-1}{1+\prior(\hat{\delta}_U-1)}(\frac{2}{3}\frac{1-\prior}{D(f_0||f_1)}+\frac{1}{3}\frac{\prior}{D(f_1||f_0)}+\bl)$-strongly convex on this interval. 
Therefore $\hat{\delta}_L$ is a unique minimizer of $C(\,\cdot\,,\hat{\delta}_U)$, and $\dl^*$ can be found by minimizing $C(\,\cdot\,,\du^*)$. $\hfill\blacksquare$


Therefore, Propositions~\ref{prop1} and~\ref{prop2} tell us we can calculate $\du^*$ explicitly as a function of $\prior$ and $\epsilon$, and $\dl^*$ numerically as the solution to a scalar, set-constrained, strongly convex optimization problem. These steps give rise to the Quickest Search Algorithm with switching Costs, which is Algorithm~1.

\begin{algorithm}[tb]
   \caption{Quickest Search Algorithm with Switching Costs}
   \label{alg1}
\begin{algorithmic}
   \STATE{\textbf{Input}: $\epsilon \in (0,1)$, $\prior \in (0,1)$, $D(f_1||f_0) > 0$, $D(f_0||f_1) > 0$}
   \STATE {$\gamma_U \gets \log \left(\frac{1-\epsilon}{\epsilon}\frac{1-\prior}{\prior}\right)$}
   \STATE{$\gamma_L \gets \log\left(\arg\min_{[0,1]}C\left(\cdot,\frac{1-\epsilon}{\epsilon}\frac{1-\prior}{\prior}\right)\right)$}
   \STATE{$t \gets 0$, $\Lambda^1_0 \gets 0$, $k \gets 1$}
   \WHILE{$\Lambda^k_t < \gamma_U$}
   \IF{$\Lambda^k_t \geq \gamma_L$}
   \STATE{\textbf{Observe}: $X^k_{t+1}$}
   \STATE{$\Lambda^k_{t+1} \gets \Lambda^k_{t} + \log\left(\frac{f_1(X^k_{t+1})}{f_0(X^k_{t+1})}\right)$}
   \STATE{$t \gets t+1$}
   \ELSE
   \STATE{$k \gets k+1$}
   \STATE{$\Lambda^k_t \gets 0$}
   \ENDIF
   \ENDWHILE
   \STATE{Label arm $k$ as a target}
\end{algorithmic}
\end{algorithm}

\section{Discussion of Brownian-Motion Approximations} \label{sec:discussion}

Now that we have presented Problem~2 as a solvable approximation of Problem~1, we will justify this substitution by showing that in the limiting cases described in Remark~1 in Section~\ref{sec:probstate}, the approximate inequalities used to formulate Problem~2 approach equalities. We are interested in the case where $f_0$ and $f_1$ are ``close'' since if they are easily distinguished, the problem is easy and optimality is not crucial.  We are interested in the case where $\epsilon$ is small because we want few errors.  Further, we are interested in the case where $\bl$ is relatively large because otherwise, the problem is solvable by the existing method of \cite{lai2011quickest}.


Recall from \cite{siegmund1985sequential} that that the inequalities \eqref{eq:gl-approx}-\eqref{eq:gl-ineq} being treated as equalities only fail to be equalities if the statistic $\Lambda^k_{t_k}$ overshoots the relevant threshold $\gamma_U$ or $\gamma_L$, rather than hitting it exactly. 
That is, while we will always have $\Lambda^k_{t_k} \geq \gamma_U$ (or $\Lambda^k_{t_k} < \gamma_L$) at the end of any stage, Problem~2 is derived by assuming $\Lambda^k_{t_k} = \gamma_U$ (or $\Lambda^k_{t_k} = \gamma_L$). This approximation is reasonable when the expected overshoot of a particular threshold is small with respect to the threshold itself, i.e., if $\mathbb{E}\left[\frac{\Lambda^k_{t_k}-\gamma_U}{\gamma_U} \mid \Lambda^k_{t_k} \geq \gamma_U \right]$ and $\mathbb{E}\left[\frac{\Lambda^k_{t_k}-\gamma_L}{\gamma_L} \mid \Lambda^k_{t_k} < \gamma_L \right]$ are small. The remainder of this section shows that 
these terms are indeed small when a problem satisfies the conditions in Remark~1.

\subsection{The case of ``close'' $f_0$ and $f_1$}

Here the phrase ``sufficiently close'' means the KL divergences $D(f_1||f_0)$ and $D(f_0||f_1)$ are small. Because $\mathbb{E}[\Lambda^k_{t_k}-\gamma_U | \Lambda^k_{t_k} \geq \gamma_U]  \leq \mathbb{E}[\Lambda^k_{t_k}-\Lambda^k_{t_k-1} | \Lambda^k_{t_k} \geq \gamma_U]$, we can see from the update law for $\Lambda^k$ in Algorithm~1 and the definition of the KL divergences that this expected overshoot approaches zero as $D(f_1 || f_0)$ and $D(f_0 || f_1)$ approach zero, as desired.
        
\subsection{The case of small $\epsilon$}

As $\epsilon$ shrinks, we must enforce a smaller error probability $P(H^{k_{\tau}} = H_0)$. From its definition, a smaller error probability directly implies a larger value of $\frac{1-\beta}{\alpha}$, which implies a larger value of $\gamma_U^*$. This relationship is intuitive: while both $\gamma_U$ and $\gamma_L$ affect $P(H^{k_{\tau}} = H_0)$, the effect of $\gamma_U$ is significantly greater since the algorithm only terminates at data stream $k$ if $\Lambda^k_{t_k} \geq \gamma_U$. 
Having a large $\gamma^*_U$ means the expected overshoots are small, as desired. 

\subsection{The case of large $\bl$}

The relationship between $\bl$ and $\gamma^*_L$ is perhaps the most interesting one. Consider the high-level goal of our analysis: to minimize the number of switches our algorithm makes before finding and identifying (hopefully correctly) a target data stream. The requirement that we find a target data stream (with probability $1-\epsilon$) means that we specifically want to avoid switching away from a target data stream. In other words, the goal is to reduce $\beta$. As with $P(H^{k_{\tau}} = H_0)$, $\beta$ depends on both $\gamma_U$ and $\gamma_L$, but the effect of $\gamma_L$ is significantly greater as the algorithm only switches away from stream $k$ if $\Lambda^k_{t_k} < \gamma_L$. Having a very negative $\gamma_L^*$ means the expected  overshoots are close to zero, as desired. 

The purpose of this analysis is to address non-trivial switching costs.
However, we do note that our rule for selecting $\gamma^*_L$ is optimal as $\bl$ approaches zero as well. Recall from Section~\ref{sec:prelim} and~\cite{lai2011quickest} that the true optimal value of $\gamma_L$ for $\bl = 0$ (i.e., the value that minimizes~\eqref{eqn:optprob}) is $\gamma_L = 0$. As $\bl \rightarrow 0$, our value of $\gamma_L^*$ found by solving Problem~2 also approaches zero, implying that the algorithm described in~\cite{lai2011quickest} is a special case of the one we develop here.

\section{Numerical Results} \label{sec:simulations}
We now compare the performance of our algorithm with the one described in~\cite{lai2011quickest}, which is the closest comparable algorithm, in a setting where switching costs are present. While the algorithm in~\cite{lai2011quickest} is optimal for the case where $\bl = 0$, it does not take into account switching costs. Furthermore, the optimal threshold $\gamma_U$ cannot be directly calculated for that algorithm, and must be estimated via Monte Carlo simulations. In contrast, our algorithm directly accounts for switching costs and uses thresholds that can be directly calculated.

In this setting, target data streams occur with prior probability $\prior = 0.1$, and obey the distribution $F_1 = \mathcal{N}(0,1)$. Nominal data streams obey $F_0 = \mathcal{N}(0,1.5)$, and we choose $\epsilon = 0.01$ to be our maximum tolerable error rate. Switching costs are drawn from a gamma distribution $\lambda_k \sim \Gamma(a,b)$, which has $\bl = \frac{a}{b}$. We will compare the performances of both algorithms across a range of values of $\bl$, by keeping $b = 1$ constant and exploring $a \in [0,5]$. The algorithm from~\cite{lai2011quickest} uses thresholds $\gamma_L = 0$ and $\gamma_U = 6.130$ for this problem regardless of the switching costs. For the algorithm described in this paper, $\gamma_U = 6.794$ regardless of switching costs, and $\gamma_L$ is chosen by solving Problem 2 with $\bl$. The values of $\gamma_L$ used and their corresponding values of $\bl$ are plotted in Figure~2.

\begin{figure}[htp]
    \centering
    \includegraphics[draft = false,width = 8cm]{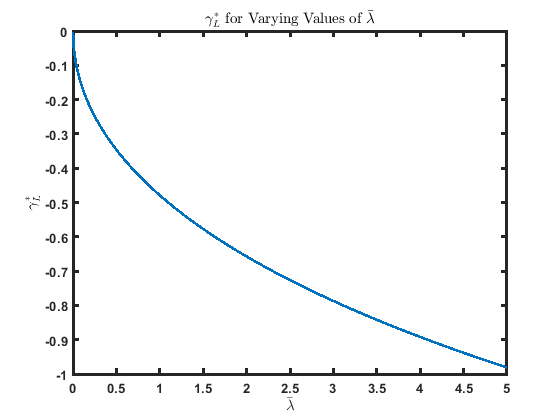}
    \caption{The larger $\overline{\lambda}$ is, the lower $\gamma_L^*$ is in Algorithm~1.}
\end{figure}

The algorithm from~\cite{lai2011quickest} achieves $\mathbb{E}_{\prior}[\tau] = 109.42$ and $\mathbb{E}_{\prior}[k_{\tau}-1] = 42.15$, and our algorithm achieves the similar numbers $\mathbb{E}_{\prior}[\tau] = 113.21$ and $\mathbb{E}_{\prior}[k_{\tau}-1] = 42.04$ for $\bl = 0$. We also note that our algorithm achieves an error rate of $0.005$. As $\bl$ grows, our goal is to reduce the combined observation/switching cost formulated in~\eqref{eqn:optprob} by reducing the number of switches. In Figure~3, we see that this is achieved. As $\gamma_L$ is varied to account for higher values of $\bl$, we see that the number of expected switches drops significantly. 

\begin{figure}[htp]
    \centering
    \includegraphics[draft = false,width = 8cm]{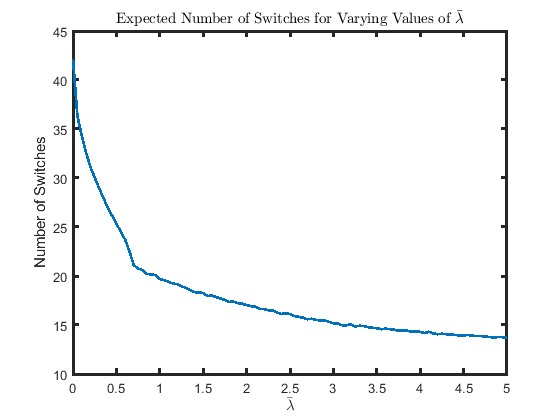}
    \caption{The expected number of switches under Algorithm~1 for different values of $\bar{\lambda}$.}
\end{figure}

The combined observation/switching costs for both algorithms are plotted in Figure~4. We can see that in the large $\bl$ region, our algorithm significantly outperforms the algorithm from~\cite{lai2011quickest}. The algorithms perform comparably up until around $\bl = 1$, after which the cost for our algorithm is always lower than the algorithm in~\cite{lai2011quickest}. Specifically, the observation/switching cost for the algorithm from~\cite{lai2011quickest} grows linearly with $\bl$. It increases by $42.15$ for every unit increase of $\bl$, since $42.15$ is the expected number of switches under that algorithm. In contrast, in the regime explored in this simulation, the use of Algorithm~1 increases the observation/switching cost at a rate of around $16.3$ per unit increase of $\bl$, meaning the cost of Algorithm~1 grows at a rate $61.3\%$ slower than the~\cite{lai2011quickest} algorithm.

\begin{figure}[htp]
    \centering
    \includegraphics[draft = false,width = 8cm]{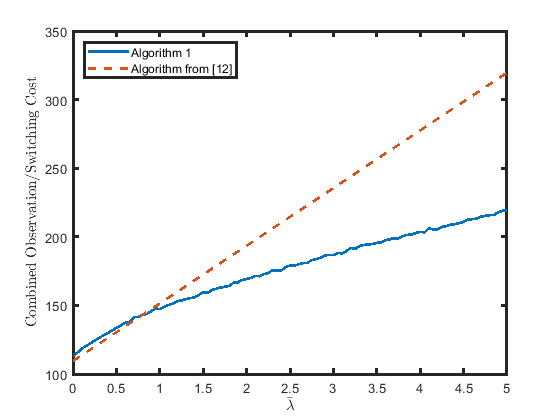}
    \caption{The combined observation/switching costs for our Algorithm~1 (blue solid line) and the algorithm from~\cite{lai2011quickest} (orange dashed line).}
\end{figure}

\section{Conclusion} \label{sec:con}
In this paper we introduced an algorithm which performs online  anomaly search over many sequences with switching costs, with parameters that can be directly calculated, and almost uniform improvement over the best comparable method that does not account for switching costs \cite{lai2011quickest}. We showed that the approximations used to derive this algorithm are accurate for problems of interest, and demonstrated the success of this algorithm with numerical simulations. Future work will embed the problem into a physical setting and perform optimal routing and control for a physical vehicle that incorporates the cost of switching that is due, e.g., to the downtime incurred by moving.

\section*{Acknowledgements}
The views and opinions expressed in this article are those of the authors and do not necessarily reflect the official policy or position of any agency of the U.S. government. Examples of analysis performed within this article are only examples. Assumptions made within the analysis are also not reflective of the position of any U.S. Government entity. The Public Affairs approval number of this document is AFRL-2023-0996.  

\bibliographystyle{IEEEtran}
\bibliography{Biblio}

\end{document}